\documentclass[11pt]{article}
\usepackage{amsmath}
\usepackage{amssymb}
\usepackage{amsfonts}
\usepackage{amsthm}

\newenvironment{Proof}{\noindent{\begin{minipage}[t]{5.9 in}}{\it Proof.  }}
{\noindent{$\diamondsuit$}{\end{minipage}}}

\theoremstyle{theorem}
\newtheorem{thm}{Theorem}[section]
\newtheorem{lemma}{Lemma}[section]
\newtheorem{prop}{Proposition}[section]
\newtheorem{conj}{Conjecture}[section]
\newtheorem{cor}{Corollary}[section]

\theoremstyle{remark}

\newtheorem{defn}{Definition}[section]
\newtheorem{remark}{Remark}[section]

\newcommand{\gp}{\mathfrak {p}}

\newcommand{\fr}{\mathop {\rm {Frob_\gp}}\nolimits}
\newcommand{\trace}{\mathop{\rm Trace}\nolimits}
\newcommand{\gal}{\mathop{\rm Gal}\nolimits}
\newcommand{\rank}{\mathop{\rm rank}\nolimits}

\newcommand{\ord}{\mathop{\rm ord}}
\newcommand{\res}{\mathop{\rm res}}
\newcommand{\Q}{{\mathbb{Q}}}
\newcommand{\C}{{\mathbb{C}}}
\newcommand{\Ql}{{\Q}_l}

\newcommand{\p}{{\mathbb{P}}}
\newcommand{\F}{{\mathbb{F}}}

\newcommand{\bk}{\bar k}
\newcommand{\pg}{\medskip\noindent}

\newcommand{\bra}[2]{\langle #1 \, , #2 \rangle_{\cH}}
\newcommand{\fra}[2]{{\frac{\scriptstyle #1}{\scriptstyle #2}}}
\newcommand{\pstar}{{\pi}^*}

\newcommand{\cS}{\cal S}

\newcommand{\cE}{\cal E}
\newcommand{\cA}{\cal A}
\newcommand{\cT}{\cal T}
\newcommand{\cF}{\cal F}
\newcommand{\cV}{\cal V}
\newcommand{\cH}{\cal H}

\newcommand{\cX}{\cal X}
\newcommand{\cY}{\cal Y}
\newcommand{\et}{{\text {\'et}}}
\newcommand{\tp}{{\tilde {\pi}}} 
\newcommand{\tS}{{\tilde {\cal S}}}

\newcommand{\tA}{{\tilde {\cal A}}}
\newcommand{\tF}{{\tilde {\cF}}}

\newcommand{\tT}{{\tilde {\cT}}}
\newcommand{\tD}{\tilde {\Delta}}
\newcommand{\tV}{{\tilde {\cV}}}
\newcommand{\kT}{{\check {\cT}}}
\newcommand{\Fp}{{\F}_\gp}

\newcommand{\bFp}{{\bar {\F}}_\gp}

\newcommand{\qp}{q_\gp}
\newcommand{\Ap}{{{\mathfrak A}_\gp}(\cA)}
\newcommand{\Bp}{{{\mathfrak B}_\gp}(\cA)}
\newcommand{\ApX}{{{\mathfrak A}_\gp}(\cX)}
\newcommand{\ApE}{{{\mathfrak A}_\gp}(\cE)}
\newcommand{\BpX}{{{\mathfrak B}_\gp}(\cX)}
\newcommand{\ap}{a_\gp}
\newcommand{\bp}{b_\gp}

\newcommand{\hK}{{\bk}({\cS})}
\newcommand{\re}{\mathop{\rm Re}\nolimits}

\newcommand{\dvs}{\mathop{\rm Div}\nolimits}

\newcommand{\pic}{\mathop{\rm Pic}\nolimits}

\newcommand{\pv}{\mathop{\rm Pic^0}\nolimits}
\newcommand{\pvs}{\mathop{\rm Pic^0_{\cS}}\nolimits}

\newcommand{\pvA}{\mathop{\rm Pic^0_A}\nolimits}
\newcommand{\pvX}{\mathop{\rm Pic^0(X/K)}\nolimits}
\newcommand{\pva}{\mathop{\rm Pic^0_{\cA}}\nolimits}

\newcommand{\pvx}{\mathop{\rm Pic^0_{\cal X}}\nolimits}
\newcommand{\ns}[1]{{\mathop{\rm NS({\it {#1}})}}}
\newcommand{\mapr}{\longrightarrow}
\newcommand{\maprlim}[1]{\smash{\mathop{\longrightarrow}\limits^{#1}}}

\newcommand{\ktr}{\mbox{${\bk}({\cS})/{\bk}$-trace}}

\title{{\bf A Local-global Summation Formula for Abelian Varieties} {\footnote{{\it keywords}: Abelian fibrations, Mordell-Weil rank, Tate conjecture, Frobenius traces}}
{\footnote{{\it MSC subject class}: 11G40, 14K15}} }
\author{\bigskip \medskip Rania Wazir {\footnote{{\it email}: wazir@dm.unito.it}} }

\begin{document}

\maketitle

\bibliographystyle{plain} 
\nocite{*}


\begin{abstract}

\pg
Let $K$ be a field finitely generated over ${\Q}$, and $A$ an Abelian variety defined over $K$.  Then by the Mordell-Weil Theorem, the set of rational points $A(K)$ is a finitely-generated Abelian group.  In this paper, assuming Tate's Conjecture on algebraic cycles, we prove a limit formula for the Mordell-Weil rank of an arbitrary family of Abelian varieties $A$ over a number field $k$; this is the Abelian fibration analogue of the Nagao formula for elliptic surfaces $E$, originally conjectured by Nagao {\cite{nag}}, and proven by Rosen and Silverman {\cite{rs}} to be equivalent to Tate's Conjecture for $E$.  We also give a short exact sequence relating the Picard Varieties of the family $A$, the parameter space, and the generic fiber, and use this to obtain an isomorphism (modulo torsion) relating the Neron-Severi group of $A$ to the Mordell-Weil group of $A$.

\end{abstract}


\section{Introduction}

\pg
At the heart of Diophantine geometry lies the link between the geometry of an algebraic variety and its arithmetic.  In this paper we explore one such link, by proving, for an Abelian fibration ${\cA}$ defined over a number field $k$, a relation between local, arithmetic information -- fibral Frobenius trace values -- and a global geometric invariant -- the Mordell-Weil rank of ${\cA}$.

\pg
To be more precise, an Abelian fibration defined over $k$ is a proper, flat morphism 
\mbox{$\pi : \cA \mapr \cS$} of smooth projective varieties, together with a section 
\mbox{$\sigma : \cS \mapr \cA$}, such that all maps and varieties are defined over $k$, and the generic fiber is an Abelian variety $A$ defined over \mbox{$K := k(\cS)$}.  We will often refer simply to ${\cA}/k$ as the Abelian fibration; when \mbox{$ \dim{\cA}=n $} and 
\mbox{$ \dim{\cS}=m $} we will also say ${\cA}/k$ is an Abelian 
\mbox{$(n,m)$-fold} defined over $k$.

\pg
If $K$ is a field finitely generated over its prime subfield, and $A$ an Abelian variety defined over $K$, then the Mordell-Weil Theorem states that the group of rational points $A(K)$ is finitely generated; in particular, this is true for $K = k(\cS)$ as above.
It is easy to show that the group of rational sections ${\cA}({\cS}/k)$ is isomorphic to the group of rational points $A(K)$, and we will call its rank the Mordell-Weil rank of $\cA$.

\pg
The Mordell-Weil rank is an important tool for studying rational points on Abelian varieties.  However, it remains a mystery, and many conjectures remain open, even in the case of elliptic curves defined over ${\Q}$.  In this paper, we focus on the Mordell-Weil rank of an Abelian fibration, and relate it to a suitable average of fibral Frobenius trace values, as follows:

\pg
For any field $F$, we denote by ${\overline F}$ its algebraic closure.  Let $O_k$ be the ring of integers of the number field $k$.  Let ${\gp}$ be a prime in $O_k$, ${\Fp}$ its residue field, and ${\qp}$ its norm (i.e., the number of elements in ${\Fp}$). For a smooth variety ${\cV}/k$, denote by ${\tV}/{\Fp}$ (or simply by ${\tV}$ where the ${\gp}$ is understood) the reduction of ${\cV}$ mod ${\gp}$.  

\pg
For $l$ prime to ${\qp}$, we consider the following Frobenius traces:
\begin{eqnarray*}
\ap{(\cV)} & := & \trace(\fr | H^1_{\et}({\cV}/{\bk}; \Ql )) \\
\bp{(\cV)} & := & \trace(\fr | H^2_{\et}({\cV}/{\bk}; \Ql )) 
\end{eqnarray*}
\noindent
This defines the Frobenius trace on non-singular fibers.  For singular fibers, set 
$$
\ap{({\cA}_x)} = \bp{({\cA}_x)} = 0.
$$
\noindent  
Furthermore, for an Abelian \mbox{$(n,m)$-fold} ${\cA}/k$, we define ''averages'' of these traces:
\begin{eqnarray*}
\Ap & := & \frac{1}{{\qp}^m} \sum_{x \in {\tS}({\Fp})} \ap{({\cA}_x)} 
= \frac{1}{{\qp}^m} \sum_{x \in {\tS}_{\rm ns}({\Fp})} \bp{({\cA}_x)} \\
\Bp & := & \frac{1}{{\qp}^m} \sum_{x \in {\tS}({\Fp})} \bp{({\cA}_x)} 
= \frac{1}{{\qp}^m} \sum_{x \in {\tS}_{\rm ns}({\Fp})} \bp{({\cA}_x)} 
\end{eqnarray*}
\noindent
where \mbox{${\tS}_{\rm ns}$} is the set of points $x$ in ${\tS}({\Fp})$ such that ${\tA}_x$ is nonsingular.

\pg
Our goal is to show:
\begin{thm}
Let 
\mbox{${\pi}: {\cA} \mapr {\cS}$} be an Abelian 
\mbox{$(n,m)$-fold} defined over $k$; assume ${\cA}(k)$ is not empty, and that
the ${{\bk}({\cS})}/{\bk}$-trace of ${A}$ is trivial.  Then Tate's Conjecture
for ${\cA}/k$ and ${\cS}/k$ implies
$$ \res_{s=1} \sum_{\gp}{-{\Ap}{\frac{\log{\qp}}{{\qp}^s}}}
+ \res_{s=2} \sum_{\gp}{{\Bp}{\frac{\log{\qp}}{{\qp}^s}}}
= \rank{\cA}({\cS}/k) + \rank \ns{A/K}.$$
\label{MThm}
\end{thm}

\pg
Recall that the \mbox{${\bk}({\cS})/{\bk}$-trace} (also called the Chow trace) of an Abelian variety $A/{\hK}$ is an Abelian variety $B/{\bk}$, together with an (injective) morphism 
\mbox{$\tau: B \mapr A$} defined over ${\hK}$ satisfying the universal mapping property for Abelian varieties $C/{\bk}$ with morphisms into $A$.  A brief summary of relevant properties is in {\cite[p.138]{sl1}}.

\begin{remark}
Although we are mainly interested in Abelian fibrations, in fact, the arguments presented in this paper can be used to prove:

\begin{thm}
Let 
\mbox{${\pi}: {\cX} \mapr {\cY}$} be a proper flat morphism of smooth projective varieties defined over $k$, with smooth generic fiber ${X}$, and section \mbox{${\sigma}: {\cY} \mapr {\cX}$}; assume ${\cX}(k)$ is not empty, and that
the ${\bk}({\cY})/{\bk}$-trace of ${\pvX}$ is trivial.  Let \mbox{$K := k({\cY})$}.  Then Tate's Conjecture
for ${\cX}/k$ and ${\cY}/k$ implies
$$ \res_{s=1} \sum_{\gp}{-{\ApX}{\frac{\log{\qp}}{{\qp}^s}}}
+ \res_{s=2} \sum_{\gp}{{\BpX}{\frac{\log{\qp}}{{\qp}^s}}}
= \rank{\pvX} + \rank \ns{X/K}.$$
\end{thm}
\end{remark}

\begin{remark}
In the case of an elliptic fibration \mbox{${\pi}: {\cE} \mapr {\cS}$} defined over $k$, with generic fiber an elliptic curve $E/K$, Theorem~{\ref{MThm}} reduces to:
\begin{equation}
\res_{s=1} \sum_{\gp}{-{\ApE}{\frac{\log{\qp}}{{\qp}^s}}}
= \rank{\cE}({\cS}/k).
\label{eqEF}
\end{equation}
\noindent
This is the analytic version of a formula conjectured by Nagao~{\cite{nag}} to hold for elliptic surfaces over ${\Q}$, and was shown by Rosen and Silverman~{\cite{rs}} to be equivalent to Tate's Conjecture for elliptic surfaces ${\cE}/k$. It was then proven in the author's thesis~{\cite{rw}} that Tate's Conjecture implies ({\ref{eqEF}}) for elliptic threefolds as well.
\end{remark}

\begin{remark}
The above equation ({\ref{eqEF}}) suggests in fact that the ${\Ap}$ sums should correspond to the Mordell-Weil rank, while the ${\Bp}$ sums correspond to the rank of Neron-Severi.  In fact, a relation of the form
\begin{equation}
\res_{s=1} \sum_{\gp}{-{\Ap}{\frac{\log{\qp}}{{\qp}^s}}}
= \rank{\cA}({\cS}/k)
\label{conj1}
\end{equation}
\noindent
can be viewed as a Birch and Swinnerton-Dyer type conjecture for families of Abelian varieties, while
\begin{equation}
\res_{s=2} \sum_{\gp}{{\Bp}{\frac{\log{\qp}}{{\qp}^s}}}
= \rank \ns{A/K}
\label{conj2}
\end{equation}
is a Tate-type conjecture, and it follows from Theorem~{\ref{MThm}} that (\ref{conj1}) and (\ref{conj2}) are equivalent.
\end{remark}

\pg
For ease of reference, we reproduce here the Birch and Swinnerton-Dyer conjecture (as generalized by Tate in 
{\cite{jt3}}), and Tate's conjecture on algebraic cycles:

\begin{conj}[Birch and Swinnerton-Dyer]
Let $A$ be an Abelian Variety defined over a global field $k$, with attached Hasse-Weil \mbox{$L$-series} $L(A/K,s)$.  Then $L(A/K,s)$ has an analytic continuation to the whole complex plane, and 
$$
\ord_{s=1} L(A/k,s) = \rank A(k)
$$
\end{conj}

\pg
\begin{conj}[Tate] {\cite[Conjecture 2]{jt1}}.
\medskip
\noindent
{\it Let ${\cV}$ be a smooth projective variety defined over $k$, and
let $L_2({\cV}, s)$ be the Hasse-Weil \mbox{$L$-function} attached to
$H^2_{\et}({\cV}/{\bk}; {\Ql})$.  Then $L_2({\cV}, s)$ has a meromorphic continuation
to ${\C}$, and has a pole at $s = 2$ of order:
$$
-\ord_{s = 2} L_2({\cV}, s) = \rank \ns{{\cV}/k}.
$$ }
\label{TC}
\end{conj}

\begin{remark}
It is also possible to consider fibrations of a surface ${\cX}$ over a curve ${\cY}$, with generic fiber a curve $X$ of genus \mbox{$g \geq 2$}.  Interesting results along these lines can be found in Wong~{\cite{sw}} and Hindry and Pacheco~{\cite{hp}}.
\end{remark}

\pg
For the remainder of this article, we fix the following notation:
\begin{defn}
Let \mbox{${\pi}: {\cA} \mapr {\cS}$} be an Abelian $(n,m)$-fold defined over $k$, with generic fiber $A$, section ${\sigma}$, and \mbox{${\Delta} := \{x \in {\cS} | {\cA}_x \; \text{is singular} \}$ } its discriminant locus.  Assume also that \mbox{${\cA}(k) \neq \emptyset$}, and that the 
\mbox{${\bk}({\cS})/{\bk}$-trace} of $A$ is trivial.  Let $K := {k}({\cS})$ be the function field of ${\cS}$, and denote by $(O)$ the image of the section $\sigma$ in ${\cA}$, and by $O$ the corresponding point on $A$.
\end{defn}

\pg
It will occasionally be necessary to throw out a finite number of "bad primes."  As we are interested in computing the residue of certain $L$-series, this will have no effect on the calculation.  Let $R$ be this set of bad primes, and start by choosing $R$ large enough such that, for all \mbox{$ {\gp} \notin R$}, \mbox{$ {\tp}: {\tA} \mapr {\tS} $} is a proper flat morphism of smooth varieties, with discriminant locus given by ${\tD}$.  All sums, products, and reductions mod ${\gp}$ will henceforth be taken only with respect to primes ${\gp} \notin R$.

\pg
We adopt the following notation for the Picard Variety, and the Picard and Neron-Severi groups of smooth varieties ${\cal X}$ defined over a field ${\kappa}$ of characteristic zero:  Let ${\pic}_{{\cal X}}$ denote the Picard Scheme of ${\cal X}/{\kappa}$; the Picard Variety of ${\cal X}/{\kappa}$, ${\pvx}$, is the component of ${\pic}_{{\cal X}}$ containing the identity, and, since \mbox{${\cal X}({\kappa}) \neq \emptyset $}, $\pvx$ is also defined over ${\kappa}$ {\cite[section 0.d]{dm2}}.  The absolute Picard and Neron-Severi groups are denoted by:
$$
\begin{array}{lcl}
{\pic}({\cal X}) & = & {{\pic}_{\cal X}}({\bar {\kappa}}) 
\text {, the Picard Group of ${\cal X}/{\bar {\kappa}}$.}\\
{\pv}({\cal X}) & = & {{\pvx}}({\bar {\kappa}}). \\
{\ns{{\cal X}}} & = & {\pic}({\cal X})/{\pv}({\cal X}) \text {, the Neron-Severi Group of 
${\cal X}/{\bar {\kappa}}$.} \\
\end{array}
$$
\noindent
Their Galois-invariant subgroups will be denoted by ${\pic}({\cX}/{\kappa})$, ${\pv}({\cX}/{\kappa})$, and ${\ns{{\cX}/{\kappa}}}$.

\medskip
\noindent
\begin{defn}
Embed $\pic(A/{\hK})$ into $\pic({\cA})$ by sending a divisor 
\mbox{$D \in \pic(A)$} to ${\bar D}$, its schematic closure in ${\cA}$.  
Fix a generating set $Z$ of $\ns{A}$,  
and let ${\bar Z}$ be its image in ${\ns{\cA}}$.
The trivial part of ${\ns{\cA}} \otimes \Q$, denoted ${\cT}$, is the 
subspace generated by ${\bar Z}$, and by all 
geometrically irreducible components of the fibral divisors.  

\pg
Denote by ${\cF}$ the subspace of ${\cT}$ generated by the
non-identity components of the fibral divisors, where the identity
component of a fibral divisor is the component intersecting $(0)$.  Then ${\cT}$ is generated by $\ns{A}$, ${\pi^*}{\ns{\cS}}$, and ${\cF}$.
\label{D5}
\end{defn}

\par
\noindent
Note that to every irreducible component of $\Delta$, there corresponds exactly one identity component of the pull-back.

\par
\noindent
\begin{remark}
For all but finitely many primes $\gp$, the trivial part ${\cT}_{\gp}$ of
\mbox{$\ns{{\tA}/{\bFp}} \otimes \, \Q$} is isomorphic to ${\tT}/{\Fp}$, the
reduction of ${\cT} \mod(\gp)$.  Therefore, by enlarging the set of bad 
primes $R$ if necessary, we can assume that 
\mbox{${\cT}_{\gp} = {\tT}/{\Fp} $}, and that \mbox{${\cF}_{\gp} = {\tF}/{\Fp}$}.  
\label{R7}
\end{remark}

\section{Galois Action on the Singular Fibers}

\newcommand{\cd}{\cal D}
\newcommand{\td}{\tilde {\cal D}}
\newcommand{\cy}{\cal Y}
\newcommand{\ty}{\tilde {\cal Y}}


\pg
Let \mbox{${\tp}: {\tA} \mapr {\tS}$} be the reduction of the fibration 
\mbox{${\pi}: {\cA} \mapr {\cS} \mod {\gp}$}.  Then \mbox{${\tp}, {\tA}, {\tS}$} are defined over ${\Fp}$, and for a point 
\mbox{$x \in {\tS}({\Fp})$}, the fiber ${\tA}_x$ is also defined over ${\Fp}$.  Let 
\mbox{$ G_{\gp} := \gal({\bFp}/{\Fp})$} be the absolute Galois group, generated by the Frobenius automorphism, ${\fr}$.  Frobenius acts on ${\tA}_x$ by sending a point 
\mbox{$\alpha \in {\tA}_x$} to the point \mbox{${\alpha}^{\qp}$}.  This has the effect of permuting the geometrically irreducible components of ${\tA}_x$, fixing only those components that are defined over ${\Fp}$ (since nonsingular fibers are irreducible, this means only singular fibers can have components that are permuted by Galois).

\pg
In this section, we find a geometric interpretation for the Galois action on the singular fibers, by showing that, on average over all fibers 
\mbox{${\tA}_x, \text{ with } x \in {\tD}({\Fp})$}, the number of fibral components fixed by Galois is equal to the number of generators of ${\tF}$ that are defined over ${\Fp}$.  Roughly, this means that every permutation possible under Frobenius is equally likely to occur for any given $x$.  We will prove this by estimating the number of rational points over the discriminant locus ${\tD}$, and showing that it is approximately \mbox{${\qp}^{(n-1)}$} times the number of generators of ${\tF}$ defined over ${\Fp}$.  As we will be using the Lang-Weil estimates {\cite{lw}}, it is convenient at this time to fix a projective embedding 
\mbox{$ {\cA} \hookrightarrow {\p}^N_k $}.

\begin{thm}  Let 
\mbox{${\pi}: {\cA} \mapr {\cS}$} be an Abelian 
\mbox{$(n,m)$-fold} defined over $k$, with notation as 
before.  Then 
$$
{\sum_{x \in {\tD}({\Fp})} \#{\tA}_x}({\Fp}) = 
{\qp}^{(n-1)}{\trace({\fr}|{\tF})} + {\qp}^{n-m}\#{\tD}({\Fp}) + O\left({\qp}^{\fra{2n-3}{2}}\right).
$$
\label{T5}
\end{thm}

\begin{Proof}
Let \mbox{${\cd} := {\pi}^{-1}({\Delta})$} be the pull-back of ${\Delta}$, with irreducible decomposition 
$$
{\cd} = {\cd}_1 + {\cd}_2 + \cdots + {\cd}_Q.
$$
\noindent
Let also ${\cy}$ be the self-intersection locus of ${\cd}$.  
Then, for all but finitely many primes $p$, ${\td}$ is the pull-back of ${\tD}$, ${\ty}$ is the self-intersection locus of ${\td}$, and thus by a Hilbert polynomial argument, dimension and degree of the varieties and their components are bounded independently of ${\gp}$.

\pg
Furthermore, observe that 
${\trace({\fr}|{\tF})}$ is the number of generators of 
${\tF}$ that are defined over ${\Fp}$, and let
\begin{eqnarray*}
J_{\gp} & := & \text{ the number of irreducible components of ${\tD}$ that are defined over ${\Fp}$,} \\
T_{\gp} & := & {\trace({\fr}|{\tF})} + J_{\gp} \leq Q.
\end{eqnarray*}
\noindent
Then $T_{\gp}$ is the number of components ${\td}_i$ that are defined over ${\Fp}$, and (by reordering if necessary) we can assume that ${\td}_i$ is defined over ${\Fp}$ for all 
\mbox{$ i \leq T_p $}.

\pg
So, determining the number of rational points on ${\td}$ reduces to:
$$
\#{\td}({\Fp}) = {\sum_{i \leq T_p} \#{\td}_i({\Fp})} + I_{\gp} +
{\sum_{i > T_p} \#{\td}_i({\Fp})},
$$
where $I_{\gp}$ is the error coming from overcounting rational points 
\mbox{$x \in {\td}_i \cap {\td}_j$}.

\pg
Now for ${i \leq T_p}$, the Lang-Weil estimates give : 
$$
\left|\#{\td}_i({\Fp}) - {\qp}^{n-1}\right| \leq C(N, n-1, g_i){\qp}^{\fra{2n-3}{2}}
$$
\noindent
for some constant $C$ depending only on $N$, $n-1$, and $g_i$, where $N$ is the dimension of the embedding space, and $(n-1)$ and $g_i$ are the dimension and degree of ${\td}_i$, respectively; as previously mentioned, these are bounded independently of 
${\gp}$, therefore 
$$
{\sum_{i \leq T_p} \#{\td}_i({\Fp})} = T_p{\qp}^{n-1} + I_{\gp} +
O\left( {\qp}^{\fra{2n-3}{2}} \right). 
$$
We use the Lang-Weil estimates again, to obtain a bound on $I_{\gp}$:
$$
I_{\gp} \leq Q{\cdot}{\#{\ty}({\Fp})} \leq Q{\cdot}C_1(N, n-2, g){\qp}^{n-2}
$$
where once again, $C_1$ is a constant depending only on $N$, $n-2$, and the genus $g$ of ${\ty}$, and therefore $I_{\gp} = O({\qp}^{n-2})$.

\pg
To see that the remaining components ${\td}_j$ for $j > T_p$ contribute only to the error term, we note that all such ${\td}_j$ are permuted by the action of ${\fr}$.  Since the rational points are fixed by ${\fr}$, they can only lie on intersections 
\mbox{${\td}_i \cap {\td}_j$}, and thus must lie in ${\ty}$.  As above, this shows that  
$$
{\sum_{j > T_p} \#{\td}_j({\Fp})} \leq  \#{\ty}({\Fp}) = O({\qp}^{n-2}).
$$
Therefore 
\begin{eqnarray*}
\#{\td}({\Fp}) & = & T_p{\qp}^{n-1} + O\left( {\qp}^{\fra{2n-3}{2}} \right) \\
 & = & {\trace({\fr}|{\tF})}{\qp}^{n-1} + J_{\gp}{\qp}^{n-1} + 
O\left( {\qp}^{\fra{2n-3}{2}} \right).
\end{eqnarray*}
By similar reasoning, we can show that
$$
\#{\tD}({\Fp})  =  J_p{\qp}^{m-1} + O\left( {\qp}^{\fra{2m-3}{2}} \right)
$$
and the theorem follows.
\end{Proof}

\begin{remark}
Although the big-O constant in Theorem~{\ref{T5}} depends on the chosen embedding, the big-O term does not contribute to the residue computation in Theorem~{\ref{MThm}}, and hence the final result is independent of choice of embedding.
\end{remark}


\section{Geometry}

\pg
We now take a closer look at the geometry of ${\cA}/k$, and extract two main theorems.  First, we will show that \mbox{${\pva}$ and ${\pvs}$} are isogenous (or equivalently, that their Tate modules are isomorphic, after Faltings).  Second, we prove a Shioda-Tate-type isomorphism, showing the relation between the Mordell-Weil rank of ${\cA}$, and the rank of the Neron-Severi group of ${\cA}$.  The results of Raynaud and Shioda in the elliptic surfaces case (see {\cite[Theorems 1 and 2]{sh1}}) were the main source of inspiration for this work. 

\pg
The basis for both theorems is an exact sequence relating the Picard Varieties of ${\cS}$, ${\cA}$, and $A$.  In order to prove this, we need to define a bilinear pairing on ${\cA}$, which will allow us to determine when a fibral divisor on ${\cA}$ is the pull-back of a divisor on ${\cS}$.  The necessary results on Intersection Theory can be found in {\cite[Appendix A]{rh}} or {\cite{wf}}.

\bigskip

\subsection{A Nondegenerate Pairing}

\pg
The intersection of two divisors \mbox{$D, F \in \dvs({\cA})$} will be an $(n-2)$-dimensional cycle.  Intersecting again with a generic \mbox{$(N+m-n+1)$-dimensional} hyperplane \mbox{${\cH} \subset {\p}^N_k $} gives an $(m-1)$-dimensional cycle on ${\cA}$, which maps to a divisor on ${\cS}$ under pushforward ${\pi}_*$.  Since the intersection pairing is symmetric and bilinear, and push-forward is linear, we have:
 
\begin{defn}
For a given \mbox{$(N+m-n+1)$-dimensional} hyperplane ${\cH}/k$ with 
\mbox{${\pi}({\cH}.{\cA}) = {\cS}$}, 
define a symmetric, bilinear pairing 
$$
\langle \:\;, \:\; \rangle_{\cH}:{\pic}({\cA}) \times {\pic}({\cA}) 
\mapr {\pic}({\cS})
$$
\noindent
via 
$$
\langle {\Lambda}, {\Upsilon} \rangle_{\cH} := {\pi}_*({\cH}.({\Lambda}.{\Upsilon}))
$$
\noindent
for any \mbox{${\Lambda}, {\Upsilon} \in {\pic}({\cA})$}.

\medskip
\par
\noindent
If $C, D$ are divisors in ${\dvs}({\cA})$, let 
\mbox{$\langle C, D \rangle_{\cH} := \langle cl(C), cl(D) \rangle_{\cH}$}, where $cl(C)$ denotes the class of $C$ modulo linear equivalence.
\label{D11}
\end{defn}
\noindent

\begin{defn}
A prime divisor \mbox{$C \subset \dvs({\cA})$} is {\bf fibral} if \mbox{$C \subset {\pi}^{-1}(G)$} for some 
\mbox{$G \subset \dvs({\cS})$}.  In general, \mbox{$D \subset \dvs({\cA})$} is a {\bf fibral divisor} if 
\mbox{$ D = {\sum_i a_iC_i}$}, where the $C_i$ are fibral.
\end{defn}

\pg
We consider first some basic facts about the pairing ${\bra{\:\;}{\:\;}}$.

\begin{enumerate}
\item
Let $G_1$, $G_2$ be distinct, irreducible divisors in ${\cS}$, and consider any 
\mbox{$D_1 \subset {\pi}^{-1}(G_1)$}, \mbox{$D_2 \subset {\pi}^{-1}(G_2)$}.  Then ${\bra{D_1}{D_2}} = 0$.  This follows because \mbox{$ D_1 \cap D_2 \subset {\pi}^{-1}(G_1 \cap G_2)$}; since $G_1$ intersects $G_2$ properly, $D_1$ intersects $D_2$ properly as well.  Intersecting again with ${\cH}$ gives a cycle of dimension $m-1$ whose image under ${\pi}$ is contained in \mbox{$G_1 \cap G_2$}, a cycle of dimension $m-2$.  By definition of ${\pi}_*$, we then have 
\mbox{${\pi}_*({\cH}.(D_1.D_2)) = 0$}.

\item
If \mbox{$cl(D) \in \pv({\cA})$}, then \mbox{${\bra{C}{D}} = 0$} for all \mbox{$C \in \dvs{\cA}$}.  This is a simple consequence of the corresponding fact about the intersection pairing.
\end{enumerate}

\begin{prop}
Let \mbox{$D \in {\dvs}({\cA})$} be a fibral divisor, and 
\mbox{$G \in {\dvs}({\cS})$}.  For ${\cH}$ any hyperplane as above, the following hold true:
\begin{itemize}
\item[(a)]
$\langle D, {{\pi}^*}(G) \rangle_{\cH} = 0$
\item[(b)]
$\langle D, D \rangle_{\cH} \leq 0$
\item[(c)]
If $\langle D, D \rangle_{\cH} = 0$ then 
$D \in {\pi}^*({\dvs}({\cS}))$.
\end{itemize} 
\label{I3}
\end{prop}

\medskip
\noindent
\begin{Proof}
\begin{itemize}
\item[(a)]
This follows from the Projection Formula, 
 once we note that for any fibral divisor $D$, \mbox{${\pi}_*({\cH}.D) = 0$}.
\begin{eqnarray*}
\langle D, {{\pi}^*}(G) \rangle_{\cH} & = &
{\pi}_*({\cH}.(cl(D).cl({{\pi}^*}{G}))) \\
  & = & {\pi}_*(({\cH}.D).{{\pi}^*}(G)) \\
  & = & {\pi}_*({\cH}.D).G  \\
  & = & 0.
\end{eqnarray*}

\medskip

\item[(b)]
Same as {\cite[Theorem 3.4b]{rw}}

\item[(c)]
Same as {\cite[Theorem 3.4c]{rw}}

\end{itemize}
\end{Proof}

\subsection{A Tate-Module Isomorphism}

\medskip
\noindent
Let \mbox{$G_k := \gal({\bk}/k)$}.  
Equipped with the bilinear pairing 
\mbox{$\langle \: \;, \: \; \rangle_{\cH}$}, we obtain:
\begin{thm}
There is an exact sequence of $G_k$-modules:
\begin{equation}
0 \mapr {\pv}({\cS}) \mapr {\pv}({\cA}) \mapr B({\bk}) \mapr 0,
\label{B6}
\end{equation}
\noindent
where $B$ is the ${\hK}/{\bk}$-trace of ${\pvA}$.
\label{SES}
\end{thm}

\par
\noindent
\begin{Proof}  

\par
\noindent
The morphism \mbox{${\pi}: {\cA} \rightarrow {\cS}$} induces
\mbox{${\pi}^*: {\pvs} \rightarrow {\pva}$} under pull-back; restriction to the generic fiber then gives
\begin{equation}
{\psi}: {{\pva} {\times}_{\cS} {\cS}} \rightarrow {\pvA}
\label{B9}
\end{equation}
\noindent
which, by the universal mapping property of the 
\mbox{${\hK}/{\bk}$-trace $(\tau,B)$} of ${\pvA}$, factors through $B$; i.e., 
there is a unique homomorphism
\mbox{$\beta: {\pva} \mapr B$} such that
$\psi = \tau \circ \beta$.
Thus, we have a sequence of morphisms

\begin{equation}
{\pvs} \maprlim{\pstar} {\pva} \maprlim{\beta} B.
\label{B2}
\end{equation}
\noindent
Considering the maps on the ${\bk}$-points, we wish to show that 
\begin{equation}
0 \mapr {\pv}({\cS}) \maprlim{\pstar} {\pv}({\cA}) \maprlim{\beta} B({\bk}) \mapr 0.
\label{B3}
\end{equation}
\noindent
is a short exact sequence of $G_k$-modules.

\pg
The injectivity of ${{\pi}^*}$ follows from the existence of the global section
\mbox{$\sigma: {\cS} \mapr {\cA}$}:  By definition, \mbox{${\pi} \circ {\sigma} = id_{\cS}$}, 
and this gives
$$id_{{\pic}({\cS})} = ({\pi} \circ {\sigma})^* = {\sigma}^* \circ {\pi}^*.$$

\pg
To show exactness at 
the middle, note that \mbox{${\psi} \circ {\pstar} = 0$}, and thus \mbox{${\beta} \circ {\pstar} = 0$} by injectivity of $\tau$.
This shows \mbox{Im($\pstar$) $\subset$ Ker($\beta$)}.

\pg
Now take any non-zero \mbox{$\Gamma \in$ Ker($\beta$)}.  
Then \mbox{$\Gamma = cl(D)$} for some
divisor $D$ with \mbox{$D|_A = 0$}, 
and thus $D$ must be a fibral divisor.  
Since \mbox{$cl(D) \in {\pv}({\cA})$}, it follows that 
\mbox{$\langle D, F \rangle_H = 0$} for every divisor $F$ on ${\cA}$,
and in particular, \mbox{$\langle D, D \rangle_H = 0$}.
By Theorem (\ref{I3}.c), this implies
\mbox{$D \in {\pi}^*({\dvs}({\cS}))$}.

\medskip
\noindent
Thus, we have $D \sim_{alg} 0$, and $D = {\pi}^*C$ for some
$C \in \dvs({\cS}).$  It remains to show that $C \sim_{alg} 0$.
But \mbox{$\sigma^*D  =  \sigma^*\pi^*C = (id_{\pic({\cS})})C =  C$}.  
Since $\sigma^*$ preserves algebraic equivalence, this shows $C \sim_{alg} 0$, and therefore 
\mbox{${\Gamma} = cl(D) \in {\pi}^*({\pv}({\cS}))$}, proving that 
\mbox{Ker($\beta$) = Im($\pstar$)}.  

\pg
Finally, to show that $\beta$ is surjective, we consider \mbox{${\Lambda} \in \pic(A {\times}_K B)$}, in the linear equivalence class of the divisor defined by 
\mbox{$\left\{ \left({\tau}(b), b \right) \in A {\times}_K B \, | \,  b \in B({\bk}) \right\}$}.  Let ${\overline {\Lambda}}$ be the schematic closure of ${\Lambda}$ in ${\cA}{\times}_k B$.  Then 
\mbox{${\overline {\Lambda}} \in \pic({\cA}{\times}_k B)$}, and by construction, there exists a dense open subset $U \subset B$ such that 
$$
{\psi}(\left. {\overline{\Lambda}} \right|_u) = {\tau}(u) \qquad \text{for all $u \in U({\bar k})$}.
$$
\noindent
Therefore, for all $u_1$, $u_2$ in $U({\bar k})$, 
$$
D_{u_1, u_2} := \left. {\overline{\Lambda}} \right|_{u_1} - \left. {\overline{\Lambda}} \right|_{u_2}
$$
is algebraically equivalent to zero, and \mbox{${\psi}(D_{u_1, u_2}) = {\tau}(u_1) - {\tau}(u_2)$}.  Fix $u_1$, and let 
$$
U_1 := \left\{ u \in B({\bar k}) \, | \, u = u_1 - u_2 \: \text{ for some $u_2 \in U({\bar k})$} \right\}.
$$  
\noindent
Then the above argument shows that ${\tau}(U_1)$ lies in the image of $pv({\cA})$ under ${\psi}$.  But since $U_1$ is open dense in $B$, it generates all of $B$, hence ${\tau}(B({\bar k})) \subset {\psi}({\pv}({\cA}))$.
\end{Proof}

\pg
It follows from the universal mapping property of the ${\bk({\cS})/{\bk}}$-trace, that isogenous Abelian varieties over ${\bk}$ have isogenous ${\bk({\cS})/{\bk}}$-traces.  Since we are assuming $A$ has trivial trace, this implies that $B$, the ${\bk({\cS})/{\bk}}$-trace of $\pva$, is also trivial.  Under these conditions, Theorem~{\ref{SES}} gives an isomorphism of 
$G_k$ modules \mbox{$\pv({\cS}) \cong \pv({\cA})$}.  Therefore, for 
all ${\ell}$ and for 
all $n$, 
$$
\pv({\cA})[{\ell}^n] \cong \pv({\cS})[{\ell}^n]
$$
and so \mbox{$T_{\ell}(\pv({\cA})) \cong T_{\ell}(\pv({\cS}))$} 
as $G_k$-modules.
\begin{thm}
$$
H^1_{\et}({\cS}/{\bk}, {\Ql}) 
\cong  H^1_{\et}({\cA}/{\bk}, {\Ql}).
$$
\noindent
as $G_k$-modules.
\label{H1ISO}
\end{thm}

\pg
\begin{Proof}
Let \mbox{$V_{\ell}(\pv({\cV})) := T_{\ell}(\pv({\cV})){\otimes}{\Ql}$}.  
Then the theorem follows from the above discussion, 
together with the Galois-invariant 
isomorphism (see~{\cite[Cor 4.19]{jm}})
$$
V_{\ell}(\pv({\cV}))(-1) \cong H^1_{\et}({\cV}/{\bk}; {\Ql})
$$
which holds for any smooth, projective variety ${\cV}$ defined over $k$.
\end{Proof}

\medskip
\noindent
\begin{cor}
${\ap}({\cA}) = {\ap}({\cS}).$
\label{CorH1ISO}
\end{cor}

\begin{remark}
In fact, from the proof it can be seen that Theorem~{\ref{SES}} is true 
for any proper, flat fibration 
\mbox{$p: {\cX} \mapr {\cY}$} of smooth, projective varieties defined over $k$, with smooth generic fiber $X$.  Furthermore, if the {\ktr} of the Picard Variety of the generic fiber $X$ is trivial, then the isomorphism in cohomology also holds, as well as the Shioda-Tate formula proven below.
\end{remark}

\subsection{A Shioda-Tate Isomorphism}

\begin{thm}[A Shioda-Tate Formula for Abelian ${\mathbf (n,m)}$-folds]

\medskip
\noindent
Embed $\pic^0(A)$ into 
$\ns{{\cA}}$ by sending a divisor $D \in \pic^0(A)$ to the divisor ${\overline D}$, where {${\overline D}$} denotes the schematic closure of $D$ in ${\cA}$.  

\medskip
\noindent
Then there is a
decomposition of $G_k$-modules,
\begin{equation}
\ns{{\cA}} \otimes \Q \cong 
\left(\pic^0(A) \otimes \Q\right) \oplus {\cT}.
\label{eqST}
\end{equation}
\label{ST}
\end{thm}

\pg
\begin{Proof}

\pg
Restriction to the generic fiber
$A$ defines a homomorphism
\begin{equation}
{\pic}({\cA}) \mapr {\pic}({A/{\hK}})
\label{B4}
\end{equation}

\noindent
which associates with every divisor class cl($D$) on ${\cA}$ the divisor 
\mbox{$D|_{A} = D.A$} on the generic fiber $A$.
Then, using the given generators for $\ns{A}$, 
\mbox{$\{Z_1, ... , Z_r\} = Z$}, write $D$ as a linear combination of $Z_i$, 
\mbox{$ D = {\sum_j n_jZ_j} $}, and adjust the image by sending 
cl($D$) to cl($D'$), where \mbox{$D' := D.E - \sum_{j=1}^r n_jZ_j$}; 
the divisor ${D'}$ is thus an element of $\pv({A})$, 
and the homomorphism becomes
\begin{equation}
{\phi}: {\pic}({\cA}) \mapr {\pv}(A/{\hK}).
\label{B5}
\end{equation}

\pg
Next, we determine the kernel of the map $\phi$:

\pg
\begin{lemma}
Let ${\kT}$ be the subgroup of ${\pic}({\cA})$ generated by the geometrically irreducible
components of the fibral divisors, and by ${\bar Z}$.  Then
\begin{equation}
0 \mapr {\kT} \maprlim{\eta} {\pic}({\cA}) \maprlim{\phi}{\pv}(A/{\hK}) \mapr 0
\label{B7}
\end{equation}

\noindent
is a short exact sequence of $G_k$-modules.
\label{L1}
\end{lemma}

\pg
{\it Proof of Lemma.}

\pg
Since \mbox{${\kT} \subset$ ker$(\phi)$} by construction, we just need to show that 
\mbox{${\kT} \supset ker(\phi)$}.  
Consider \mbox{${\Upsilon} \in$ ker$(\phi)$}, 
i.e., \mbox{${\Upsilon} = cl(D)$}, 
where \mbox{$D|_A - \sum_{j=1}^r n_jS_j \sim 0$ on $A$}.
But then
$$
D|_A- \sum_{j=1}^r n_jS_j  = div(h),
$$
\noindent
where \mbox{$ h \in {\hK}(A) = {\bk}({\cA})$}, and hence there exists \mbox{$H \in {\bk}({\cA})$} such that 
\mbox{$(H)|_A = (h)$}.

\pg
If \mbox{$D' := D - \sum_{j=1}^r n_j{\bar Z}_j - (H)$}, 
then $D'$ must be in some fiber, i.e.,
\mbox{$D' \in {\kT}$}, and therefore, since also \mbox{${\bar Z} \in {\kT}$},
$$
{\Upsilon} = cl(D) = cl(D') + {\bar Z} \in {\kT}.
$$

\pg
Finally, given any 
\mbox{${\hK}$-rational} divisor
$C$ on $A$, taking the schematic closure of its irreducible
components gives a divisor ${\bar C}$ on ${\cA}$ such that
\mbox{${\bar C}.A = C$}.  Therefore ${\phi}$ is surjective. 

\pg
Under the embedding \mbox{$\pvA \mapr \ns{\cA}$}, $\pvA$ inherits a $G_k$ action, such that restriction to the generic fiber is $G_k$-invariant.  Therefore also the map $\phi$ is $G_k$-invariant.  ${\diamondsuit}$

\pg
To complete the proof of Theorem~{\ref{ST}}, note that the exact sequence~({\ref {B2}}) shows that restriction to the generic fiber sends $\pva$ to $B = 0$.  
Since \mbox{$ \ns{{\cA}} := {\pic}({\cA})/{\pva} $}, combining this with the short exact sequence~({\ref {B7}}) and tensoring with ${\Q}$ gives the desired isomorphism.
\end{Proof}

\pg
Taking \mbox{$ \gal ({\bk}/k)$-invariants} of (\ref {eqST}), and noting that $A$ and ${\pvA}$ are isogenous Abelian Varieties, hence have the same rank, we get as corollary the following 
formula relating $\rank\ns{{\cA}/k}$ and $\rank{\cA}({\cS}/k)$:
\begin{cor}
$$
\rank\ns{{\cA}/k} = \rank{\cA}({\cS}/k) + \rank\ns{{\cS}/k} + \rank\ns{A/K}
+ \rank{\cF}^{G_k}.
$$
\label{CorST}
\end{cor}


\section{The Main Theorem}

\newcommand{\errm}{O\left({\qp}^{\fra{2(n-m)-3}{2}}\right)}
\newcommand{\errn}{O\left({\qp}^{\fra{2n-3}{2}}\right)}
\newcommand{\tr}{{\qp}^{(n-1)}{\trace({\fr}|{\tF})}}
\newcommand{\trn}{{\trace({\fr}|{\tF})}}

\par
\noindent
Before proceeding to the proof of Theorem~{\ref{MThm}}, we review briefly the definitions of the $L$-series we will be dealing with. 
Recall that, for a smooth variety ${\cV}/k$, 
the {\bf Hasse-Weil ${\mathbf L}$-series} attached to 
$H^2_{\et}({\cV}/{\bk}; {\Ql})$, denoted $L_2({\cV}, s)$, is given by
$$
L_2({\cV}, s) := 
\prod_{\gp} \det\left(1 - {\fr}{\qp}^{-s} | H^2_{\et}({\cV}/{\bk}; {\Ql})\right)^{-1}.
$$
\medskip
\noindent
Let $V$ be a finite-dimensional ${\Q}$-vector space, with an action of $G_k := \gal({\bk}/k)$.  Then $V$ defines a Galois representation of $G_k$, and the {\bf Artin ${\mathbf L}$-series} attached to $V$ is {\cite[VII.10]{jn}}:
$$
L({V}, s) := 
\prod_{\gp} \det\left(1 - {\fr}{\qp}^{-s} | V^{G_k} \right)^{-1}.
$$
\begin{remark}
To be precise, since in this paper we are working over all primes
${\gp} \notin R$,  
$$
L_2({\cE}, s)  \approx  \prod_{{\gp} \notin R} \det\left(1 - {\fr}{\qp}^{-s} | H^2_{\et}({\cE}/{\bk}; {\Ql})\right)^{-1}, 
$$
\noindent
where the symbol $\approx$ is used to indicate that the two sides
agree up to finitely many Euler factors.  Similarly for $L_2({\cS}, s)$ and $L({\cF}, s)$.  
This, however, has no
effect on the residue computation.
\end{remark}

\pg
To determine the precise relation between the order of vanishing of these $L$-series, and the Frobenius traces, we now take logarithmic derivatives, and compute residues:

\pg
Consider first $L({\cF}, s)$, with \mbox{$\re(s) > \frac{1}{2}$}
\begin{eqnarray*}
\frac{{\rm d}}{{\rm {ds}}} \log{L({\cF}, s)} & = & 
\frac{{\rm d}}{{\rm {ds}}} \sum_{\gp} -\log 
\det\left(1 - {\fr}{\qp}^{-s} | {\cF}\right) \\
 & = & \sum_{\gp} -\trace\left(\fr | {\cF}\right)\frac{\log{\qp}}{{\qp}^s}
+ O(1) .
\end{eqnarray*}
\noindent
Therefore,
\begin{eqnarray}
\res_{s=1} \sum_{\gp} \trace({\fr}|{\tF}){\frac {\log{\qp}}{{\qp}^s}}
& = & -\res_{s = 1}\frac{{\rm d}}{{\rm {ds}}} \log{L({\cF}, s)} \nonumber \\ 
& = & - \ord_{s=1} L({\cF}, s) \nonumber \\
 & = & \rank({\cF}^{G_k}),
\label{E11}
\end{eqnarray}
\noindent
where this last equality follows from {\cite[Proposition 1.5.1]{rs}}. 

\pg
Furthermore, for $L_2({\cA}, s)$ with $\re(s) > \frac{3}{2}$,
\begin{eqnarray*}
\frac{{\rm d}}{{\rm {ds}}} \log{L_2({\cA}, s)} & = & 
\frac{{\rm d}}{{\rm {ds}}} \sum_{\gp} -\log 
\det\left(1 - {\fr}{\qp}^{-s} | H^2_{\et}({\cA}/{\bk}, \Ql)\right) \\
 & = & \sum_{\gp} -{\bp}({\cA})\frac{\log{\qp}}{{\qp}^s}
+ O(1) .
\end{eqnarray*}
\noindent
Therefore,
\begin{eqnarray}
\res_{s=2} \sum_{\gp} {\bp}({\cA}){\frac {\log{\qp}}{{\qp}^s}}
& = & -\res_{s = 1}\frac{{\rm d}}{{\rm {ds}}} \log{L_2({\cA}, s)} \nonumber\\  
& = & -\ord_{s=2} L_2({\cA}, s) \nonumber \\
 & = & \rank \ns{{\cA}/k} \quad \text {by Tate's Conjecture},
\label{E13}
\end{eqnarray}
\noindent
and similarly
\begin{eqnarray}
\res_{s=2} \sum_{\gp} {\bp}({\cS}){\frac {\log{\qp}}{{\qp}^s}}
& = & -\res_{s = 1}\frac{{\rm d}}{{\rm {ds}}} \log{L_2({\cS}, s)} \nonumber\\  
& = & -\ord_{s=2} L_2({\cS}, s) \nonumber \\
 & = & \rank \ns{{\cS}/k} \quad \text {by Tate's Conjecture}.
\label{E12}
\end{eqnarray}

\medskip
\noindent
Finally, to prove Theorem~{\ref{MThm}}, we will use a counting argument to determine (up to an error term of order ${\qp}^{\left({\fra{2n-3}{2}}\right)}$ ), an equation for ${\Ap}$ and ${\Bp}$ in terms of ${\ap}$ and ${\bp}$ of ${\cA}$ and ${\cS}$, plus a term coming from the singular fibers.  Pulling together the results of the previous sections, combined with the $L$-series interpretations above, will complete the proof.

\setcounter{section}{1}
\setcounter{thm}{0}
\medskip
\noindent 
\begin{thm}
Let 
\mbox{${\pi}: {\cA} \mapr {\cS}$} be an Abelian 
\mbox{$(n,m)$-fold} defined over $k$; assume ${\cA}(k)$ is not empty, and that
the ${\hK}/{\bk}$-trace of ${A}$ is trivial.  Then Tate's Conjecture
for ${\cA}/k$ and ${\cS}/k$ implies
$$ \res_{s=1} \sum_{\gp}{-{\Ap}{\frac{\log{\qp}}{{\qp}^s}}}
+ \res_{s=2} \sum_{\gp}{{\Bp}{\frac{\log{\qp}}{{\qp}^s}}}
= \rank{\cA}({\cS}/k) + \rank \ns{A/K}.$$
\end{thm}
\setcounter{section}{6}

\begin{Proof}

\pg
By the Lefschetz Fixed-Point Theorem, the number of rational points on 
${\tA}$ and ${\tS}$ is given by:
\begin{eqnarray}
\#{\tA}({\Fp}) & = & O\left({\qp}^{\fra{2n-3}{2}}\right) 
+{\qp}^{(n-2)}{\bp}({\cA})
-{\qp}^{(n-1)}{\ap}({\cA})+{\qp}^n.
\label{E2} \\
\#{\tS}({\Fp}) & = & O\left({\qp}^{\fra{2m-3}{2}}\right) 
+{\qp}^{(m-2)}{\bp}({\cS})-{\qp}^{(m-1)}{\ap}({\cS})+{\qp}^m,
\label{E3}
\end{eqnarray}

\pg
It is also possible to determine the number of rational points on ${\tA}$ by counting fiber by fiber: 
\begin{eqnarray}
\#{\tA}({\Fp}) & = & {\sum_{x \in {\tS}({\Fp})} \#{\tA}_x({\Fp})} \nonumber \\
& = & {\sum_{x \in {\tS}_{ns}({\Fp})} \#{\tA}_x({\Fp})} + {\sum_{x \in {\tD}({\Fp})} \#{\tA}_x({\Fp})}.
\label{eqC1}
\end{eqnarray}

\pg
By Theorem~{\ref{T5}}, we have 
\begin{equation}
{\sum_{x \in {\tD}({\Fp})} {\tA}_x} = 
{\qp}^{(n-1)}{\trace({\fr}|{\tF})} + {\qp}^{n-m}\#{\tD}({\Fp}) + O\left({\qp}^{\fra{2n-3}{2}}\right).
\label{eqC3}
\end{equation}

\pg
To determine 
\mbox{${\sum_{x \in {\tS}_{ns}({\Fp})} \#{\tA}_x({\Fp})}$}, we note that 
${\tA}_x$ is smooth for all \mbox{$x \in {\tS}_{ns}({\Fp})$}, hence the Lefschetz Fixed Point Theorem applies again:
$$
\#{\tA}_x({\Fp}) = {\qp}^{(n-m)} - {\ap}({\cA}_x){\qp}^{(n-m-1)} + {\bp}({\cA}_x){\qp}^{(n-m-2)} + {\errm}
$$
Therefore
\begin{eqnarray}
{\sum_{x \in {\tS}_{ns}({\Fp})} \#{\tA}_x({\Fp})} & = & 
{\qp}^{(n-m)}{\#{\tS}_{ns}({\Fp})} - {\Ap}{\qp}^{n-1} + {\Bp}{\qp}^{n-2} + {\errn} \nonumber \\
& = & {\qp}^{(n-m)}({\#{\tS}({\Fp})} - {\#{\tD}({\Fp})} ) - {\Ap}{\qp}^{n-1} + {\Bp}{\qp}^{n-2} 
+ {\errn} \nonumber\\
& = & {\qp}^n - {\ap}({\cS}){\qp}^{(n-1)} + {\bp}({\cS}){\qp}^{(n-2)} 
- {\#{\tD}({\Fp})}{\qp}^{(n-m)}\nonumber \\
&   &  - {\Ap}{\qp}^{n-1} + {\Bp}{\qp}^{n-2} + {\errn} 
\label{eqC2}
\end{eqnarray}
where this last equality is obtained by substituting the expression in (\ref{E3}) for 
\mbox{$\#{\tS}({\Fp})$}.  Combining with (\ref{eqC3}), we have an alternative expression for ${\#{\tA}({\Fp})}$:
\begin{eqnarray}
\#{\tA}({\Fp}) & = & {\qp}^n - {\ap}({\cS}){\qp}^{(n-1)} + {\bp}({\cS}){\qp}^{(n-2)} - {\Ap}{\qp}^{n-1} + {\Bp}{\qp}^{n-2} \nonumber \\
& & + {\tr} + {\errn} \label{eqC4}
\end{eqnarray}
Finally, equating the two expressions ({\ref{E2}}) and ({\ref{eqC4}}) for the number of rational points on ${\tA}$, we obtain:
\begin{eqnarray}
& & - {\Ap}{\qp}^{n-1} + {\Bp}{\qp}^{n-2} \\
& = & {\ap}({\cS}){\qp}^{(n-1)} - {\ap}({\cA}){\qp}^{(n-1)} 
+{\bp}({\cA}){\qp}^{(n-2)} - {\bp}({\cS}){\qp}^{(n-2)} \nonumber \\
& & \quad - {\tr} + {\errn} \nonumber \\
& = & {\bp}({\cA}){\qp}^{(n-2)} - {\bp}({\cS}){\qp}^{(n-2)} - {\tr} + {\errn}, 
\qquad \text{by Corollary~{\ref{CorH1ISO}}}. \nonumber
\label{E6}
\end{eqnarray}

\pg
Summing over ${\gp}$ and taking residues, we have
\begin{eqnarray}
& & \res_{s=1} \sum_{\gp} -{\Ap}{\frac {\log{\qp}}{{\qp}^s}} + 
\res_{s=2} \sum_{\gp} {\Bp}{\frac {\log{\qp}}{{\qp}^s}} \\
& = & \res_{s=2} \sum_{\gp} {\bp}({\cA}){\frac {\log{\qp}}{{\qp}^s}}
 - \res_{s=2} \sum_{\gp} {\bp}({\cS}){\frac {\log{\qp}}{{\qp}^s}} 
- \res_{s=1} \sum_{\gp} {\trn}{\frac {\log{\qp}}{{\qp}^s}}   \nonumber \\
 & = & - \ord_{s=2} L_2({\cA}, s) + \ord_{s=2} L_2({\cS}, s) + \ord_{s=1} L({\cF}, s) 
\label{eqORD}
\end{eqnarray}
\noindent
where the final equality comes from (\ref{E11}), (\ref{E13}), and (\ref{E12}).  Applying Tate's Conjecture to (\ref{eqORD}) above leaves:
$$
\res_{s=1} \sum_{\gp} -{\Ap}{\frac {\log{\qp}}{{\qp}^s}} + 
\res_{s=2} \sum_{\gp} {\Bp}{\frac {\log{\qp}}{{\qp}^s}}
 =  - \rank({\cF}^{\gal({\bk}/k)}) - \rank \ns{{\cS}/k} 
+ \rank \ns{{\cA}/k}.
$$
\noindent
By the Shioda-Tate formula for Abelian fibrations 
(Corollary {\ref {CorST}}), this implies
\begin{equation}
\res_{s=1} \sum_{\gp} -{\Ap}{\frac {\log{\qp}}{{\qp}^s}}
 + \res_{s=2} \sum_{\gp} {\Bp}{\frac {\log{\qp}}{{\qp}^s}}
= \rank {\cA}({\cS}/k) + \rank \ns{A/K}.
\label{eqFIN}
\end{equation}

\end{Proof}

\begin{cor}
Let \mbox{${\pi}: {\cA} \mapr {\cS}$} satisfy the hypotheses of Theorem~{\ref{MThm}}, and assume in addition that $L_2({\cA}, s)$ and $L_2({\cS}, s)$ do not vanish on the line 
\mbox{$\re{(s)} = 2$}.  Then
$$
\lim_{X {\mapr} \infty} {\frac{1}{X}}\sum_{{\qp} \leq X} -{\Ap}{\log{\qp}}
 + \lim_{X {\mapr} \infty} {\frac{1}{X}}\sum_{{\qp} \leq X} {\Bp}{\frac{\log{\qp}}{\qp}}
= \rank {\cA}({\cS}/k) + \rank \ns{A/K}.
$$
\end{cor}

\begin{Proof}

\pg
By {\cite[Proposition 1.5.1]{rs}}, $L({\cF}, s)$ does not vanish on the real line 
\mbox{$\re{(s)} = 1$}.  Furthermore, the non-vanishing assumption on $L_2({\cA}, s)$ and  $L_2({\cS}, s)$ means we have a Tauberian theorem for these $L$-series, such that
\begin{eqnarray}
\lim_{X {\mapr} \infty} {\frac{1}{X}}\sum_{{\qp} \leq X} -{\bp}({\cA}){\frac{\log{\qp}}{\qp}}
& = & \ord_{s=2} L_2({\cA}, s) \\
\lim_{X {\mapr} \infty} {\frac{1}{X}}\sum_{{\qp} \leq X} -{\bp}({\cS}){\frac{\log{\qp}}{\qp}}
& = & \ord_{s=2} L_2({\cS}, s) \\
\lim_{X {\mapr} \infty} {\frac{1}{X}}\sum_{{\qp} \leq X} -{\trn}{\log{\qp}}
& = & \ord_{s=1} L({\cF}, s)
\end{eqnarray}
Therefore
$$
\lim_{X {\mapr} \infty} {\frac{1}{X}}\sum_{{\qp} \leq X} -{\Ap}{\log{\qp}}
 + \lim_{X {\mapr} \infty} {\frac{1}{X}}\sum_{{\qp} \leq X} {\Bp}{\frac{\log{\qp}}{\qp}}
= - \ord_{s=2} L_2({\cA}, s) + \ord_{s=2} L_2({\cS}, s) + \ord_{s=1} L({\cF}, s)
$$
and the result follows by applying Tate's Conjecture and the Shioda-Tate formula as before.  \end{Proof}

\clearpage

\bibliography{main} 
 
\end{document}